\documentstyle[12pt]{amsart}
\def\merto{\hbox{$-\!\rightarrow$}}
\newcommand{\PP}{{\mathbb P}}
\newcommand{\R}{{\mathbb R}}
\newcommand{\wt}{\widetilde}
\renewcommand{\tilde}{\widetilde}
\newcommand{\wh}{\widehat}
\newcommand{\epq}{{\cal E}^{p,q}}
\newcommand{\cA}{{\cal A}}
\newcommand{\cC}{{\cal C}}
\newcommand{\cD}{{\cal D}}
\newcommand{\cE}{{\cal E}}
\newcommand{\cL}{{\cal L}}
\newcommand{\cM}{{\cal M}}
\newcommand{\cO}{{\cal O}}
\newcommand{\cP}{{\cal P}}

\newcommand{\loc}{{\operatorname{loc}}}
\newcommand{\codim}{{\operatorname{codim}}}

\newcommand{\Supp}{{\operatorname{Supp}}}
\newcommand{\rank}{{\operatorname{rank}}}
\newcommand{\Id}{{\operatorname{Id}}}
\newcommand{\Image}{{\operatorname{Image}}}
\newcommand{\dee}{\partial}
\newcommand{\deebar}{\overline{\partial}}
\newcommand{\ddc}{dd^c}
\renewcommand{\epsilon}{\varepsilon}
\renewcommand{\phi}{\varphi}
\newcommand{\eqd}{\buildrel {\operatorname{def}}\over =}
\newcommand{\nm}{\min(n,m)}

\newcommand{\C}{{\mathbb C}}

\newcommand{\Z}{{\mathbb Z}}
\newcommand{\G}{{\mathbb G}}
\newcommand{\glm}{\G(\ell,m)}

\title
[Value Distribution for Sequences and Complex
Dynamics]
{Value Distribution for Sequences of Rational Mappings and Complex
Dynamics}

\author{Alexander Russakovskii}

\address{\hskip-\parindent
Alexander Russakovskii\\
Theory of Functions Department\\
Institute for Low Temperature Physics\\
310164 Kharkov\\ 
Ukraine}
\email{russakov@@math.jhu.edu, russakov@@msri.org}

\author{Bernard Shiffman}

\address{\hskip-\parindent
Bernard Shiffman\\
Department of Mathematics\\
Johns Hopkins University \\
Baltimore, MD 21218 \\
USA}
\email{shiffman@@math.jhu.edu, shiffman@@msri.org}

\thanks{Shiffman's research  was
partially supported by National Science Foundation Grant DMS-9500491.
The research of both authors at MSRI is
partially supported by National Science Foundation Grant 
DMS-9022140.}
\date{March 28, 1996}

\curraddr{Mathematical Sciences Research Institute, 
1000 Centennial Drive, Berkeley, CA 94720, USA (both authors, until
June 1996)}

\newtheorem{theorem}{{\sc Theorem}}[section]
\newtheorem{cor}[theorem]{{\sc Corollary}}
\newtheorem{lemma}[theorem]{{\sc Lemma}}
\newtheorem{proposition}[theorem]{{\sc Proposition}}

\newenvironment{proof}{\noindent{\em Proof:\/}}{}
\newenvironment{example}{\smallskip\noindent{\it
Example:\/}}{\smallskip}
\newenvironment{rem}{\smallskip\noindent{\it Remark:\/}}{\smallskip}
\newenvironment{defin}{\smallskip\noindent{\it Definition:\/}}{\smallskip}
\begin{document}

\thispagestyle{empty}

\begin{abstract} We obtain results on the asymptotic equidistribution of the
pre-images of linear subspaces for sequences of rational mappings between
projective spaces.  As an application to complex dynamics, we consider the
iterates $P_k$ of a rational mapping $P$ of
$\PP^n$. We show, assuming a condition on the topological degree
$\lambda$ of
$P$, that there is a probability measure
$\mu$ on
$\PP^n$ such that the discrete measures $\lambda^{-k}P_k^*\delta_w$ converge
to $\mu$ for all
$w\in\PP^n$ outside a pluripolar set. 
\end{abstract} 
\maketitle

\section{Introduction}

The study of value distribution for sequences of mappings may be regarded as
an analogue of Nevanlinna theory; instead of studying the asymptotic behavior
of the area of a pre-image of an analytic set in a ball when its radius tends
to infinity, one investigates the asymptotics of the pre-images of an analytic
set under a sequence of mappings.  One of the main reasons for studying this
subject is its applications to complex dynamics, although value distribution
theory for sequences is also of independent interest.

Many investigations have the Brolin-Lyubich Theorem (see
\cite{Brolin,Lyubich,FreireLopesMane}) as their starting point. This theorem
can be formulated as follows:

\medskip

{\it Let $R(z)$ be a rational function of degree  $d \geq 2,$
and let $R_k$ be its $k$-th iterate.
Then there is an invariant probability measure $\mu$ on
$\PP^1=\C\cup\{\infty\}$ such that for all $w$ outside an exceptional
set $\cE\subset\PP^1$ containing at most $2$ points,

$$\frac{1}{d^k} (R_k)^* \delta_w \rightarrow \mu.$$
}

Here $(R_k)^* \delta_w$ is a discrete measure in $\PP^1$, counting
the number of roots of the equation $R_k (z) = w$ with
multiplicities; convergence is in the weak sense.

Many of the recent papers on holomorphic dynamics in several complex variables
(see \cite{BedfordSmillie1,BedfordSmillie2,BedfordSmillie3,
BLS1,BLS2,HubbardPapadopol,FornaessSibony,FS0,FS1,FS2}) have dealt with
various extensions of this theorem to $n$ variables in particular cases.  This
result was also extended by Sodin \cite{Sodin} to the case where $\{R_n\}$ is
an arbitrary sequence of rational functions of one variable with rapidly
increasing degrees.  Sodin proved that if one agrees to omit a larger set
$\cE$ of exceptional values, then the pre-images of the remaining values
$w\in\PP^1\setminus\cE$ are equidistributed in a certain sense.

In a recent paper \cite {RS}, the phenomenon of equidistribution for sequences
of polynomial mappings $\C^n\to \PP^m$ was studied.  The phenomenon and the
corresponding results may be described as follows. We let $\omega=\omega _n$
denote the K\"ahler form of the Fubini-Study metric on complex projective
$n$-space $\PP^n$ (normalized so that $\int_{\PP^{n}}\omega
^n=\nobreak 1)$.

\medskip
\begin{it}
Suppose that $\{P_k\}$ is a sequence of polynomial mappings $\C^n
\rightarrow \PP^m$, such that $\sum \frac{1}{\delta_1(P_k)}
< \infty$, where $\delta_1(P_k)$ denotes the maximal degree of the components
of $P_k$.  Then
\begin{enumerate}
\item[(i)] pre-images of all but an exceptional pluripolar set of  complex
hyperplanes in $\PP^m$  are equidistributed with the pull-backs $P_k^*
\omega$.

\item[(ii)] if all $P_k$ are non-degenerate, then pre-images of all but a
pluripolar set of points in $\PP^m$ are equidistributed with $P_k^*
\omega^m \ (m \leq n).$
\end{enumerate} 
\end{it}

The exceptional sets in (i) and (ii) above are described in terms of a
``proximity sequence.''

In \cite{RS}, only the cases of codimension $1$ and $m$ were treated. The
present paper deals with all intermediate codimensions $1\leq\ell\leq \nm$ and
rational mappings $\PP^n\to \PP^m$ (Theorem~\ref{th1}) as well as giving a
higher dimensional version (with pluripolar exceptional set) of the
Brolin-Lyubich theorem (Theorem~\ref{th2}).

In order to state our results on the asymptotic equidistribution of
pre-images of linear subspaces of intermediate dimension, we let $\G(\ell,
m)$ denote the Grassmannian of projective linear subspaces of codimension
$\ell$ in $\PP^m$. Note that $\G(m,m) = \PP^m$. If $P:\PP^n\merto\PP^m$ is a
meromorphic (i.e., rational) map, we let $\delta_\ell(P)$ denote the degree of
$P^{-1}(W)$ for generic $W\in\glm$ ($1\leq\ell\leq \nm$). (For generic $W$,
$P^{-1}(W)$ has pure codimension $\ell$.)  One easily sees that $\delta_1(P)$
is the degree of the polynomials in a representation of $P$ (using homogeneous
coordinates); if $P$ is holomorphic (i.e., regular), then
$\delta_\ell(P)=\delta_1 (P)^\ell$.  In Section~\ref{s-degrees}, we give an
analytic description of the ``intermediate degrees'' $\delta_\ell(P)$ and we
show that $\delta_{k+\ell}(P)\leq\delta_k(P)\delta_\ell(P)$. Intermediate 
degrees have been considered also in \cite{Friedland}.

Our first result on the equidistribution of pre-images under an arbitrary
sequence of rational maps of projective spaces is as follows:

\bigskip 
\begin{theorem}  Let $\{P_k\}$ be a sequence of rational mappings
from $\PP^n$ to $\PP^m$.  Let $1\leq\ell \leq \nm$, and let $\{a_k\}$ be a
sequence of positive numbers  such that
$$\sum_{k=1}^\infty\frac{\delta_{\ell-1}  (P_k)}{a_k}<+\infty\;.$$  Then
there exists a pluripolar subset $\cE$ of
$\glm$ such that
$$\frac{1}{a_k}\left(P^*_k
[W]-P^*_k\omega^\ell\right)\rightarrow 0$$ as $k\rightarrow \infty$
for all $W\in\glm\setminus\cE$.
\label{th1}
\end{theorem}

Here, convergence means weak convergence in $\cD'{}^{\ell,\ell}(\PP^n)$. In
Section~\ref{s-exceptional}, we give a description of the exceptional 
set $\cE$ in terms of the
``proximity sequence.'' The pull-back $P^*_k\omega^\ell$ is smooth off the
indeterminacy locus of $P$ and has locally integrable coefficients.
For generic $W$,
$P^*_k[W]$ is given by integration over $P_k^{-1}(W)$. Precise definitions of
the pull-back currents $P^*_k[W], P^*_k\omega^\ell\in\cD'{}^{\ell,\ell}(\PP^n)$
are given in Section~\ref{s-notation}.

As an application of Theorem~\ref{th1} to complex dynamics, we
have the following result on the equidistribution of iterated pre-images:

\begin{cor}  Let $P: \PP^n \merto \PP^n$ be a rational mapping,
and let $P_k$ denote the $k$-th iterate of $P$.  If
$a > \delta_{\ell-1}(P)$, then
$$\frac{1}{a ^k} \left(P^*_k [W] - P^*_k
\omega ^\ell\right)\rightarrow 0$$ for all $W\in\glm$ outside a pluripolar set.
\label{cor1}
\end{cor} 

\noindent To obtain Corollary~\ref{cor1} from Theorem~\ref{th1}, we use the
fact that $\delta_{\ell-1}(P_k) \leq \delta_{\ell-1}(P)^k$, which we verify in
Section~\ref{s-degrees} (see Lemma~\ref{qp}).

For holomorphic maps $P:\PP^n \to \PP^n$ of degree $d\geq 2$ it was 
shown in \cite {HubbardPapadopol} and \cite{FS2} that $\frac{1}{d^k}\log
|P_k|^2$ converges uniformly to a Green's function $G$ (where $P_k$ is
the $k$-th iterate of $P$), and hence by Bedford-Taylor \cite{BT},
$$\frac{1}{d^{nk}}P^*_k\omega ^n\rightarrow(dd^c G)^n= \mu\;,$$
where $\mu$ is a probability measure on $\PP^n$ which is invariant in the
sense that $P_* \mu = \mu$.

However, Fornaess and Sibony \cite{FS2} showed that for a non-holomorphic
map $P:\PP^2 \merto \PP^2$ of degree $d$, the number of points 
in the generic fibre of $P$ is strictly less than $d^2$. (We give a
generalization of this fact as Lemma~\ref{algebra}.)  Thus, for such $P$
it follows that $$\frac{1}{d^{2k}}P^*_k\omega ^2\rightarrow 0\;.$$ So the
question of existence of a nontrivial limit measure in the meromorphic case
remained open.  The following result gives a limit measure and
equidistribution for the iterates of meromorphic maps on $\PP^n$:

\begin{theorem} Let $P:\PP^n \merto \PP^n$ be a
rational map, let $\lambda=\delta_n(P)$ denote the topological degree of $P$,
and write  
$$\mu_k=\frac{1}{\lambda^k}P_k^*\omega^n\;.$$
If $\lambda> \delta_{n-1}(P)$,  then the sequence $\{\mu_k\}$ converges to
a probability measure $\mu$ on $\PP^n$ and
$$\frac{1}{\lambda ^k} P^*_k [W]\rightarrow \mu$$
for all points $W\in \PP^n$ outside a pluripolar
set.
\label{th2}\end{theorem}  

\begin{rem} We have for any test function $\phi$,
$$(P_*\mu_{k+1}, \phi) = \frac{1}{\lambda^{k+1}} \int 
(\phi\circ P) P_{k+1}^*\omega^n
=\frac{1}{\lambda^{k+1}} \int P^*
(\phi P_{k}^*\omega^n) 
=\frac{1}{\lambda^{k}} \int 
\phi P_{k}^*\omega^n = (\mu_k, \phi)\;.$$

Thus \begin{equation} P_* \mu_{k+1}=\mu_k\;.\label{last}\end{equation}
If $\mu(I_P)=0$,
then $P_*\mu$ would be
well-defined and hence we could let $k\to\infty$ in (\ref{last}) to conclude
that $\mu$ is an invariant measure, i.e., $P_*\mu=\mu$.
\end{rem}

Further results will be given in a subsequent paper.

%%%%%%%%%%%%%%%%%%
\section {Examples and open questions\label{s-examples}}

\bigskip

We begin with several elementary examples of  sequences of iterated
mappings illustrating our results. These examples are the extensions to
$\PP^2$ of proper polynomial mappings of $\C^2$.
Besides equidistribution, we are going to look at invariance
properties of the limit currents and measures. 

\noindent{\it Example 1.\/} (This example was discussed in \cite {RS}.)
Let $P_k (z,w)$ be the $k$-th iterate of the mapping
$P=(z^\delta,w^\delta):\C^2 \rightarrow \C^2$. We first restrict our attention
to $\C^2$. It is easy to see that the uniform convergence
$\displaystyle\frac{1}{\delta^k} \log (1+|P_k (z,w)|^2) \rightarrow G(z,w)$
takes place, where $$G(z,w) = \sup (\log^+ |z|^2, \log^+ |w|^2)$$ is the
plurisubharmonic Green function of the unit polydisk with logarithmic growth
at infinity. It follows that
$$\ddc G = \lim_{k\to\infty}\frac{1}{\delta^k} P_k^* \omega.$$
The current $T=\ddc G$ is concentrated on the boundary of the
unit polydisk and on the set $\{ |z| = |w|>1\}.$
This current is invariant under the mapping $P$ in the sense that
$$P^* T = \delta T\;.$$ (The pull-backs of closed positive $(1,1)$-currents
by holomorphic and meromorphic maps are defined in Section~\ref{s-notation}.)
According to the results of \cite{RS}, this current is the limit of
pre-images of all nonexceptional hyperplanes (lines). 

Every hyperplane of the form $z = c$ or $w = c$ has pre-images of similar
form. These pre-images tend towards the cylinder $|z| = 1$ or
$|w|=1$, respectively. Thus all such hyperplanes are exceptional. A third
family of exceptional hyperplanes consists of those passing through the
origin. Pre-images of these hyperplanes tend towards the ``cone" $\{|z| =
|w|\}$.  We now consider
the complex projective plane $\PP^2=\C^2\cup H_\infty$, identifying
$(z,w)\in\C^2$ with $(1,z,w)\in\PP^2$. We let
$\PP^{2*}\approx\PP^2$ denote the parameter space of hyperplanes in
$\PP^2$; the point $(\zeta_0,\zeta_1,\zeta_2)\in\PP^{2*}$ represents the
hyperplane $\{(z,w)\in\C^2: \zeta_0+ \zeta_1 z + \zeta_2 w=0\}$. (The point
$(1,0,0)$ represents the hyperplane at infinity,
$H_\infty$.)  Thus we see that the set of exceptional hyperplanes consists of
the three pencils in $\PP^{2*}: \ \{\zeta_0 = 0\},\ \{\zeta_1 = 0\},\
\{\zeta_2 = 0\}$. We note that, besides $T$, there are at least five more
(linearly independent) invariant closed $(1,1)$-currents (on $\C^2$),
$T_1=\ddc \log^+ |z|^2,\ T_2=\ddc \log^+|w|^2, \ T_3 = \ddc \log \max (|z|^2,
|w|^2),\ T_4=[H_1]=\ddc\log |z|^2,\ \ T_5=[H_2]=\ddc\log |w|^2$ (where
$H_1=\{z=0\},H_2=\{w=0\}$) with the same property
$$P^* T_i = \delta T_i.$$  In fact, they are limits of pre-images of the
corresponding exceptional hyperplanes. (These currents are invariant on
$\PP^2$ as well as on
$\C^2$, and we have the additional invariant current
$T_6=[H_\infty]$ on  $\PP^2$.)

Now consider pre-images of points. The measure $\mu=T^2=(\ddc G)^2,$ being the
limit of $\displaystyle\frac{1}{\delta ^{2k}} P_k^* \omega^2,$ is concentrated
on the distinguished boundary of the unit polydisk; hence pre-images of most
points must tend to the torus by Theorem~2 of \cite {RS}. However, all points
of the form $(0, c)$ or $(c, 0)$ have pre-images of the form $(0, c')$
(respectively $(c', 0)$) and are definitely exceptional (as well as the points
of $H_\infty$). So in this case the exceptional set is the
same union of 3 hyperplanes in $\PP^{2*}.$

The measure $\mu$ possesses invariance properties, 
$P^* \mu = \delta ^2\mu.$ Note that we also have $\mu = T_1 \wedge T_2.$

The map $P$ extends to a rational mapping: $Q(t,z,w)=(t^\delta , z^\delta ,
w^\delta ):\PP^2\to\PP^2$, which we call the {\it projectivization} of $P$.
(For this example, $Q$ is holomorphic.) The situation then becomes very
symmetric with respect to all variables. To simplify the terminology, we will
call the corresponding currents (measures) ``projectivizations'' of the ones
defined on $\C^n.$ Instead of, say, $\ddc \log^+|z|^2=\ddc \log (|z|^2 \vee
1)$ one has to consider $\ddc \log (|z|^2 \vee |t|^2),$ where $\vee$ stands
for maximum, and so on.  The same points and hyperplanes are exceptional, and
because of the absence of the indeterminacy set, no difficulty occurs.

\medskip
\noindent
{\it Example 2.\/} The situation changes if we consider the iterates of the
mapping $P=(z^{d_1},w^{d_2}): \C^2 \rightarrow \C^2$ with $d_1\neq d_2.$ For
simplicity, let $P=(z^2, w^3).$ It is easy to see that in this case
$\displaystyle\frac{1}{3^k} \log (1+|P_k (z,w)|^2) \rightarrow G(z,w)$, where
now $G(z,w) = \log^+ |w|^2.$ It follows as before that
$$\ddc G = \lim_{k\to\infty}\frac{1}{3^k} P_k^* \omega;$$ however this time
the current $T=\ddc G$ is concentrated on the cylinder $\{ |w|=1\}.$ It is
invariant for the mapping $P$: $$P^* T = 3T\;.$$ Applying the results of
\cite{RS}, we see that this current is the limit of the pre-images of all
nonexceptional hyperplanes.

The set of exceptional hyperplanes in $\C^2$ consists of the two pencils in
$\PP^{2*}: \ \{\zeta_0 = 0\},\ \{\zeta_2 = 0\}$.  Besides $T$, there are at
least three more invariant currents: $T_1=\ddc \log^+ |z|^2$, $T_2 = \ddc
\log |z|^2$,  $T_3 = \ddc\log |w|^2$. However this time
$$P^* T_1 = 2T_1\;.$$ (Note that for any smooth form
$S$ on $\C^2$, we have $P_*P^*S=6S$; by smoothing and taking limits, one
sees that this identity is also valid if $S$ is a positive
$(1,1)$-current on $\C^2$. Hence
$P_*T_1=3T_1,\ P_*T=2T$.)  The currents $T_i$ are the limits of
pre-images of the corresponding families of exceptional hyperplanes. The
coefficient 2 for
$T_1$ makes the situation somewhat different as we shall see.

Now consider pre-images of points. The results of \cite{RS} do not provide
useful information since now $T^2=(\ddc G)^2\equiv 0.$ According to
Theorem~\ref{th2}, this is not a surprise, since $\delta_2(P) = 6$, not 9, and
we have a limit
$$\mu=\lim_{k\to\infty}\frac{1}{6^k} P_k^* \omega^2,$$ which is the
same measure (concentrated on the distinguished boundary of the unit
polydisk) as in Example~1. So pre-images of most points must concentrate
there.  Again, the exceptional set for points is the same union of 3
hyperplanes in
$\PP^2.$
The measure $\mu$ possesses the invariance property 
$P^* \mu = 6\mu$. Note that we also have $\mu = T \wedge T_1.$

Consider the projectivization of this mapping: $Q(t,z,w)=(t^3, t z^2,
w^3):\PP^2\to\PP^2$ and the sequence $Q_k$ of its iterations.  Note that
although $P$ is holomorphic, $Q$ is only meromorphic and has one indeterminacy
point $(0,1,0)$ which is the reason for all ``anomalies.''  Since $H_\infty$
is contracted by $Q$ to the fixed point $(0,0,1)$, the map $Q$ is ``generic''
in the sense of Fornaess and Sibony \cite{FS1,FS2}.  Also, the graph $G$ of
the mapping $Q$ is singular. So it is really necessary to resolve the
singularities of $G$ to define all our currents correctly (see
Section~\ref{s-notation}). Note that we have $\delta_1(Q_k)= 3^k,\
\delta_2(Q_k)=6^k<\delta_1^2(Q_k)$. In accordance with the results of
\cite{FS1,FS2}, there is an invariant $(1,1)$-current $T=\ddc \log (|w|^2 \vee
|t|^2)$ which is the projectivization of the above mentioned current on
$\C^2.$ The projectivization of the other invariant current, $T_1$, is no
longer invariant, since now $Q^*T_1 = 2T_1 + H_\infty$.

The situation is even more complicated for the point case. The results of
\cite {FS1,FS2} do not provide a nontrivial invariant measure since $T^2=0.$
According to Theorem~\ref{th2}, there is a limit measure $\mu'$ in this case.
In fact, $\mu'=[\ddc \log (|z|^2 \vee |w|^2 \vee
|t|^2)]^2$ is the projectivization of $\mu$ and is invariant on $\PP^2$. We
have $\mu'=T
\wedge T_1$ as before, although now $T_1$ is not invariant itself.  The
measure $\mu'$ has the property $Q^*\mu' = 6 \mu'$.

\medskip\noindent
{\it Example 3.\/} Finally, consider the sequence of iterations of
$P(z,w)=(w^3, z^2).$ This example is not generic in the sense of Fornaess and
Sibony
\cite{FS1,FS2}, since $H_\infty$ is contracted by (the projectivization of)
$P$ to the indeterminacy point
$(0,1,0)$, so their results cannot be applied. Since $\delta_1 (P_k)
\neq \delta_1^k(P)$ (which is the case for non-generic maps) we cannot
conclude that the sequence of currents $\frac{1}{\delta_1(P_k)} P_k^* \omega$
converges. Instead we have two subsequences
$$\{P_{2k} = (z^{6^k}, w^{6^k})\}\;,\quad \{ P_{2k+1} = (w^{3\cdot 6^k},
z^{2\cdot 6^k})\}\;.$$ The first subsequence is the same as in Example~1 with
$\delta = 6$, and the corresponding current subsequence thus converges to
$$T_1=\ddc \log^+ (|z|^2 \vee |w|^2).$$ The second subsequence of currents
converges to $$T_2=\ddc
\log^+ (|w|^2 
\vee |z|^{4/3})\;.$$
These currents are responsible for the distribution of pre-images of 
hyperplanes. As for the invariance properties, there seems to be no
invariant current since we have $P^* T_1 = 3 T_2,\ P^* T_2 = 2 T_1.$

For the codimension two case, we have $\delta_2 (P_k) = \delta_2^k(P)
=6^k$. So there is a limit measure which is the same as in the two previous
examples and is invariant.

\medskip We state here some open problems.
In the Brolin-Lyubich Theorem, the exceptional set consists of two
points. In the above examples, the exceptional sets
are unions of at most $3$ hyperplanes in $\PP^2$. This leads to the questions:

{\it Can we further describe the exceptional set in
Theorem~\ref{th2}?  Is the exceptional set in fact algebraic?}

Another question is:  

{\it Does the measure $\mu$ in Theorem~\ref{th2} charge the indeterminacy set
of $P$?}

If the answer to this question is ``no'', then by the remark following Theorem~\ref{th2}, $P_*\mu$ would be well-defined
and $\mu$ would be an invariant measure.

%%%%%%%%%%%%%%%%%%
\section {Notation and terminology \label{s-notation}}

We let $\epq(X)$, $\cD^{p,q}(X)$, $\cD'{}^{p,q}(X)$ denote the
spaces of (complex-valued)
$\cC^\infty$ forms, compactly supported $\cC^\infty$ forms, and currents,
respectively, of bidegree
$(p,q)$ on a complex manifold $X$ and we use the standard differentials
$d=\dee +\deebar$,
$d^c=(4\pi \sqrt{-1})^{-1}(\dee - \deebar )$. Points in complex
projective $n$-space $\PP^n$ are identified with their
representations $z = (z_0, z_1,\ldots, z_n)$ in homogeneous coordinates. 
We shall regard the Grassmannian $\glm$
of projective linear subspaces of codimension $\ell$
in $\PP^m$
as a subvariety of 
$\PP(\textstyle{\bigwedge^{m+1-\ell}}\C^{m+1})$.
 
If $P:\PP^n \merto \PP^m$ is a non-constant meromorphic map, it is a
well-known consequence of Chow's theorem that $P$ must be rational, i.e., $P$
can be written in the form $P = (P_0,\ldots, P_n)$ where $P_j\in
\C[z_0,\ldots,z_n]$ and $\deg P_0=\cdots=\deg P_m=d$.  We can assume that the
$P_j$ have no common factors; we then say that $\deg P=d$. (This notion of
degree should not be confused with the {\it topological degree} of an
equidimensional rational map, which is the number of points in
the pre-image of a generic point in the range.) We let
$I_P\subset\PP^n$ denote the indeterminacy locus of $P$ (the points where
$P$ is not holomorphic); $I_P$ is an algebraic subvariety of codimension
$\geq 2$.

Suppose $f:Y\rightarrow X$ is a holomorphic mapping of complex manifolds.  If
$\alpha$  is a current on
$X$, $f^*\alpha$  is not always defined.  However, we shall define $f^*\alpha
$ in two special cases:  First, we suppose $\alpha =u\gamma$ where $\gamma
\in \cE^{p,q}(X)$ is a smooth form and $u$ is the difference of
plurisubharmonic functions.  Assume further that $f(Y)$ is not contained in
the $\pm
\infty$ locus of $u$.  Then $u\circ f$ is the difference of plurisubharmonic
functions on $Y$ and hence is in ${\cal L}^1_{\loc}(Y)$.  We define
$f^*\alpha =(u\circ f)f^*\gamma$, which is clearly independent of the
representation $\alpha =u\gamma$.  The second case we consider is that of a
current of the form
$[D]\wedge \gamma$, where $[D]$ is the current of integration over a divisor
$D$ on $X$ and
$\gamma$ is a smooth form as before.  We assume also that $f(Y)\not\subset
\Supp\: D$ so that
$f^*D$ is a divisor on $Y$.  We then define
$f^*([D]\wedge\gamma) = [f^*D]\wedge f^*\gamma$.  These two definitions are
consistent in the following way.  Suppose $\alpha =\log |g|^2\cdot\gamma$
where $g$ is a meromorphic function on $X$ such that neither the zeroes nor
the poles of $g$ contain $f(Y)$ and $\gamma$ is a closed $(p,q)$-form on $X$. 
Then
$dd^c\alpha =[D]\wedge\gamma$ where $D=$ Div$(g)$.  Hence
\begin{equation} f^*dd^c \alpha = [f^*D]\wedge f^*\gamma=dd^c(\log|g\circ f|^2 \cdot
f^*\gamma) = dd^cf^*\alpha \;.\label{trivial}\end{equation}

Let $P:\PP^n\merto\PP^m$ be a rational map.  For a smooth $(p,q)$-form $\eta
\in
\epq(\PP^m)$ we define the pull-back current $P^*\eta\in
\cD'{}^{p,q}(\PP^m)$ as follows:
We let $G_P\subset \PP^n\times \PP^m$ denote the graph of $P$ (which is an
irreducible algebraic subvariety of $\PP^n\times\PP^m$) and we consider a
desingularization $\wt G \stackrel{\rho}{\rightarrow }G_P$.  We have the
commutative diagram:
\begin{equation}\begin{array}{lcr}
&\wt G&\\[.1cm]\; {\pi_1}\swarrow & &\!\searrow{\pi_2}\\[.1cm]
\PP^n&\stackrel{P}{\merto}&\PP^m
\end{array}\label{comm-diagram}\end{equation}
We then define
$$P^*\eta = \pi_{1*} \pi^*_2 \eta\;.$$ Note that $P^*\eta$ has coefficients
in ${\cal L}^1_{\loc}$ and has singular support contained in the indeterminacy
locus $I_P$ of $P$.  In fact,
$(P|_{\PP^{n}\setminus I_{P}})^*\eta$ is the usual pull-back of the form
$\eta$.

We consider the current of integration $[W]\in\cD'{}^{\ell,\ell}(\PP^m)$ and
define the pull-back $\pi_2^*[W]$ to be the current of integration over the
algebraic $(n-\ell)$-cycle $\pi_2^*W$ on $\wt G$ (using the diagram
(\ref{comm-diagram})), whenever $\dim \pi_2^{-1}(W)=n-\ell$.  If we represent
$W$ as the intersection of hyperplanes $H_1,\ldots,H_{\ell}$ in $\PP^m$, then
$\pi_2^*W$ is the intersection of divisors $\pi_2^*H_1\cap\cdots\cap\pi_2^*
H_\ell$.  For a definition of this intersection, which is a formal sum of the
irreducible components of $\pi_2^{-1}(W)$ with positive integer coefficients,
see \cite[Appendix~A]{Ha} or Definitions 2.3 and 2.4.2 (or Example 7.1.10) in
\cite{Fu}.  This pull-back, or intersection, can also be defined analytically
as follows. Let $g_j$ be a local defining function for $\pi_2^*H_j$, for
$1\leq j\leq \ell$, and write $g=(g_1,\ldots,g_\ell)$.  Then by Griffiths and
King \cite[1.10]{GK} (see also \cite[I.12, Th.~3]{Shabat}), we have the local
formula $$\pi_2^*[W]=\ddc\left(\log |g|^2(\ddc\log |g|^2)^{\ell-1}\right)\;.$$
(Alternately, $\pi_2^*[W]=dd^c\log|g_1|^2\wedge\ldots\wedge
dd^c\log|g_\ell|^2$, where the existence of this product of currents is
guaranteed by Demailly \cite{De}.)  One way to verify that these definitions
are all equivalent is to first consider generic hyperplanes
$H_1,\ldots,H_{\ell}$ so that the divisors $\pi_2^* H_1,\ldots, \pi_2^*
H_{\ell}$ are smooth hypersurfaces (of multiplicity 1) in $\wt {G}$
intersecting transversely. (This is possible by Lemma~\ref{Bertini} in the next
section.)  Then for this case, $\pi_2^*[W]$ is the current of integration over
a smooth submanifold.  In the general case, the current $\pi_2^*[W]$ is the
weak limit of the pull-backs of generic intersections
$W^\nu=H_1^\nu\cap\cdots\cap H_\ell^\nu$ converging to $W$.  The existence of
the limit follows, for example, from \cite[Cor.~11.1]{Fu} for the algebraic
definition and from \cite{De} or \cite[Cor.~3.6]{FS5} for the analytic
definition.

We now state the Poincar\'e-Lelong formula for linear subspaces of
$\PP^m$ and describe its pull-backs by a rational map $P:\PP^n\merto\PP^m$. Let
$W\in\glm$ be an $(m-\ell)$-plane in $\PP^m$. For each $W\in\glm$, we define
the current
$$\Lambda_W = \log \frac{|\zeta|^2|W|^2}{|\zeta\wedge W|^2}\sum^{\ell-1}_{j=0}
(d_\zeta d^c_\zeta \log|\zeta\wedge W|^2)^j\wedge \omega _\zeta ^{\ell-1-j}
\in \cD'{}^{\ell-1, \ell-1}(\PP^m_\zeta)\;,$$ which (by definition) has locally
integrable coefficients.  We have the {\it generalized Poincar\'e-Lelong
formula} for $W$ \cite[1.15]{GK} (see also \cite[II.6, pp.~68-69]{Shabat}),
\begin{equation} dd^c\Lambda_W=\omega^\ell-[W]\;.\label{Poincare-Lelong}\end{equation} Now let $W =
H_1\cap\cdots\cap H_{\ell}$ such that $\dim \pi_2^{-1}(W)=n-\ell$, where we
use the notation of (\ref{comm-diagram}).  Applying the generalized
Poincar\'e-Lelong formula (\cite[1.15]{GK} or \cite[II.6, pp.~68-69]{Shabat}),
to the divisors $\pi_2^* H_1,\ldots, \pi_2^* H_{\ell}$ of the lifted
hyperplane-section bundle $\pi_2^*\cO_{\PP^m}(1)$ with Chern form
$\pi_2^*\omega$, we obtain
\begin {equation}
\ddc \pi_2^* \Lambda_W = \pi_2^* \omega^{\ell} - \pi_2^*[W]\;,
\label{lifted-Poincare-Lelong}\end{equation} where
$\pi_2^*\Lambda_W\in\cD'{}^{\ell-1,\ell-1}(\wt G)$ is given by
\begin{equation}(\pi_2^* \Lambda_W)(\wt z) \eqd \log \frac{|\pi_2(\tilde{z})|^2
|W|^2}{|\pi_2(\tilde{z}) \wedge W|^2} \sum\limits_{j=0}^{\ell -1} (\ddc \log
|\pi_2(\tilde{z}) \wedge W|^2)^j \wedge \pi_2^* \omega^{\ell -1 - j} \in
\cD'{}^{\ell-1,\ell-1} (\wt{G}).\label{lambda1}\end{equation}

In particular, $\pi_2^* \Lambda_W$ has $\cL^1_{\loc}$ coefficients and 
is smooth on $\wt{G}\setminus \pi_2^{-1}(W)$. We
define the currents
\begin{equation}\begin{array}{rcl}P^*\Lambda_W &=&\pi_{1*}(\pi_2^* \Lambda_W)\;,\\
P^*[W]&=&\pi_{1*}[\pi_2^{-1} W]\;.\end{array}\label{lambda2}\end{equation}
By applying $\pi_{1*}$ to (\ref{lifted-Poincare-Lelong}), we 
obtain the current identity on $\PP^n$,
\begin{equation}
\ddc P^*\Lambda_W = P^*\omega^{\ell} - P^*[W].
\label{usual-Poincare-Lelong}\end{equation}
Note that for generic $W$, $\pi_2^*[W]$ has multiplicity identically 1 and
contains no components inside the exceptional locus of $\pi_1$, and thus
$[P^*W]$ is the current of integration over the closure of
$(P|_{\PP^n\setminus I_P})^{-1}(W)$.

%%%%%%%%%%%%%%
\section{The intermediate degrees of a rational map\label{s-degrees}}

In this section, we give some properties of the intermediate degrees
$\delta_\ell(P)$ of a rational map $P$, which we also describe analytically
and topologically. We use the following consequence of Bertini's theorem:

\begin{lemma}
Let $Y$ be a projective algebraic manifold and 
let $f:Y\to\PP^m$ be a nonconstant holomorphic map.  Then for a generic
hyperplane $H\subset\PP^m$, the divisor $f^*H$ is smooth and has
multiplicity $1$.
\label{Bertini}\end{lemma}

\begin{proof} Apply Bertini's theorem (see, for example,
\cite[p.~137]{GH}) to the complete linear system of $f^*H$.
\hfill $\Box$\end{proof}

For a subvariety $V\subset\PP^m$, we write $P^{-1}(V)=\pi_1(\pi_2^{-1}(V))$
(using the notation in (\ref{comm-diagram})).  We let $\#(S)$ denote the
cardinality of a set $S$. We begin with a formula for the integral of certain
singular forms on $\PP^n$:
\begin{lemma}
Let $P_j:\PP^n \merto \PP^{m_{j}}$, $1\leq 
j\leq n$, be rational maps, and let  $I=I_{P_1}\cup\cdots\cup I_{P_n}$.
\begin{itemize}
\item[\rm i)] For generic hyperplanes
$H_1\subset\PP^{m_1},\ldots,H_n\subset\PP^{m_n}$,
$$\int_{\PP^n\setminus I} P^*_1 \omega _{m_{1}}\wedge\cdots\wedge P^*_n \omega
_{m_{n}}=\#\left(\bigcap_{j=1}^n P_j^{-1}(H_j) \setminus I\right)\;,$$
\item[\rm ii)]
$$\int_{\PP^n\setminus I} P^*_1 \omega _{m_{1}}\wedge\cdots\wedge
P^*_n \omega _{m_{n}}\leq \prod^n_{j=1} \deg
P_j\;.$$\end{itemize}
\label{general formula}\end{lemma}

To verify Lemma~\ref{general
formula}, we first give a topological description of the integral.  Let
$X=\PP^{m_1}\times\cdots\times\PP^{m_n}$ and write
$\cP=(P_1,\ldots,P_n):\PP^n\merto X$; then $I=I_\cP$.  Consider the
commutative diagram
$$\begin{array}{lcr}
&\wt G&\\[.1cm]\; {\pi_1}\swarrow & &\!\searrow{\pi_2}\\[.1cm]
\PP^n&\stackrel{\cP}{\merto}&X
\end{array}$$
where $\wt G$ is a desingularization of the graph of $\cP$ and
$\pi_1,\pi_2$ are the projections. Let $p_j:X\to\PP^{m_j}$ denote the
projection to the $j$-th factor and let $\wt P_j=p_j\circ\pi_2:\wt
G\to\PP^{m_j}$, for
$1\leq j\leq n$.  Let
$t_m$ denote the positive generator of $H^2(\PP^m,\Z)$, and write
$$t_{m_1}\times\cdots\times t_{m_n}=p_1^*t_{m_1}\bullet\cdots\bullet
p_n^* t_{m_n}\in H^{2n}(X,\Z)\;,$$ where $\bullet$ denotes the cup
product in the cohomology ring.

\begin{lemma}  Using the notation of Lemma~\ref{general formula},
\begin{eqnarray*}\displaystyle\int_{\PP^n\setminus I} P^*_1 \omega
_{m_{1}}\wedge\cdots\wedge P^*_n \omega
_{m_{n}}&=&\left(\pi^*_2(t_{m_1}\times\cdots\times t_{m_n}),\wt
G\right)\\&=&(\wt P_1^*t_{m_1}\bullet\cdots\bullet\wt
P_n^*t_{m_n},\wt G)\in\Z\;.\end{eqnarray*}
\label{topology}\end{lemma}

\begin{proof}  Since the de Rham class of the K\"ahler form $\omega_m$
on $\PP^m$ equals $t_m$, we have by Section~\ref{s-notation},
\begin{eqnarray*}\int_{\PP^n\setminus I} P^*_1 \omega
_{m_{1}}\wedge\cdots\wedge P^*_n \omega
_{m_{n}}&=&\int_{\PP^n\setminus
I_\cP}\cP^*(\omega_{m_1}\times\cdots\times\omega_{m_n})\\
&=&\displaystyle\left(\pi_{1*}\pi^*_2(\omega_{m_1}
\times\cdots\times\omega_{m_n}),1\right)
\\&=&\displaystyle\left(\pi^*_2(\omega_{m_1}
\times\cdots\times\omega_{m_n}),
1\right)\\&=&\left(\pi^*_2(t_{m_1}\times\cdots\times
t_{m_n}),\wt G\right)\;.\end{eqnarray*}
The second equality follows from our definitions.\hfill $\Box$\end{proof}

\noindent{\it Proof of Lemma~\ref{general formula}:\/} Let $E\subset\wt G$ be
the exceptional locus of $\pi_1$.  By Lemma~\ref{Bertini} applied to $\wt
P_1:\wt G\to\PP^{m_1}$, there is a hyperplane $H_1\subset\PP^{m_1}$ such that
the divisor $\pi^*_2(H_1\times\PP^{m_2}\times\cdots\times\PP^{m_n})$ is a
smooth hypersurface $Y_1\subset\wt G$ of multiplicity 1 with $\dim Y_1\cap
E<n-1$.  By Lemma~\ref{Bertini}, we can inductively find hyperplanes
$H_2\subset\PP^{m_2},\ldots, H_n\subset\PP^{m_n}$ such that, writing
$$Y_j=\pi_2^{-1}(H_1\times\cdots\times
H_j\times\PP^{m_{j+1}}\times\cdots\times\PP^{m_n})\;,$$ $Y_j$ is a smooth
submanifold of $\wt G$ of dimension $n-j$, $\dim
Y_j\cap E<n-j$, and the divisor $(\wt P_j|_{Y_{j-1}})^*H_j$ (on $Y_{j-1}$)
has multiplicity 1, or equivalently,
$\wt P_j^{-1}(H_j)$ intersects $Y_{j-1}$ transversely.
In particular, $\dim Y_n=0$ and $Y_n\cap E=\emptyset$.
Write $H'_j=p^{-1}_j(H_j)$, for $1\leq j\leq n$, so that
$$Y_j=\pi_2^{-1}(H'_1\cap\cdots\cap H'_j)=\pi_2^{-1}(H'_j)\cap Y_{j-1}\;.$$  

A codimension $j$ submanifold $S$ of a complex manifold $Y$ determines 
the current of integration $[S]\in\cD'{}^{j,j}(Y)$ given by
$([S],\phi)=\int_S\phi$ for a test form $\phi$.  We also let $[S]$
denote the de Rham class in $H^{2j}(Y,\R)$ containing the current
$[S]$.  (If $Y$ is compact of dimension $n$, then the cohomology class
$[S]$ is the Poincar\'e dual of the $(2n-2j)$-cycle $S$.)  If two
submanifolds
$S_1,S_2$ intersect transversely, then $[S_1\cap S_2]=[S_1]\bullet [S_2]$
in the cohomology ring of $Y$.
In particular, $[H_j]=t_j\in H^2(\PP^{m_j},\Z)\subset H^2(\PP^{m_j},\R)$
and
$$[H_1\times\cdots\times H_n]=[H'_1\cap\cdots\cap
H'_n]=t_{m_1}\times\cdots\times t_{m_n}\;.$$  Furthermore, our
construction of the $H_j,Y_j$ implies that
$$[Y_j]= [Y_{j-1}]\bullet\pi_2^*[H'_j]=\pi_2^*[H'_1]\bullet\cdots\bullet
\pi_2^*[H'_j]=\pi^*_2[H'_1\cap\cdots\cap H'_j]\;.$$
Therefore,
\begin{equation}\pi^*_2(t_{m_1}\times\cdots\times
t_{m_n})=\pi^*_2[H_1\times\cdots\times
H_n]=[Y_n]\label{a}\end{equation}
(where the points of $Y_n$ have multiplicity 1).

Let $A_j:\C^{m_j+1}\to\C$ be a linear map
defining the hyperplane $H_j$, and consider the polynomial
$$Q_j=A_j(P_{j0},\ldots,P_{jn})\in\C[z_0,\ldots,z_n]\;,$$ where
$P_j=(P_{j0},\ldots,P_{jn})$, for $1\leq j\leq n$.  We then have
\begin{equation}\{z\in\PP^n:Q_j(z)=0\}=\pi_1\left(\pi_2^{-1}(H'_j)
\right)\supset I_{P_j}\;.\label{b}\end{equation}  Since $Y_j\subset\wt G\setminus E$
and
$\pi_1$ maps $\wt G\setminus E$ bijectively to $\PP^n\setminus I_\cP$, we
have
\begin{eqnarray*}\pi_1(Y_n)&=&\pi_1\left(\pi_2^{-1}
(H'_1\cap\cdots\cap H'_n)\right)\setminus I_\cP\\&=&\bigcap_{j=1}^n
\pi_1\left(\pi_2^{-1}(H'_j)\right)\setminus I_\cP=\bigcap_{j=1}^n
P_j^{-1}(H_j)\setminus I_\cP\\ &=&\{z\in\PP^n\setminus
I_\cP:Q_1(z)=\cdots=Q_n(z)=0\}\;.\end{eqnarray*} Thus by (\ref{a}),
$$\left(\pi^*_2(t_{m_1}\times\cdots\times t_{m_n}),\wt
G\right)=\left([Y_n],\wt G\right)=\#(Y_n)\;.$$  By B\'ezout's Theorem,
$$\#(Y_n)=\#(\pi_1(Y_n))=\#\left(\bigcap_{j=1}^n P_j^{-1}(H_j)\setminus I_\cP
\right)\leq\prod_{j=1}^n\deg Q_j =\prod_{j=1}^n\deg P_j\;.$$
The conclusion follows from Lemma~\ref{topology}. \qed

\begin{defin} Let $P:\PP^n\merto\PP^m$ be a rational map.  We define the {\it
intermediate degrees}  $\delta_\ell (P)$ of $P$ by the formula
$$\delta_\ell(P) =
\int_{\PP^n\setminus I_P}P^*\omega_m^\ell \wedge
\omega _n^{n-\ell}$$ for $1\leq\ell \leq \nm$.
\end{defin}

We shall show in Lemma~\ref{algebra} below that the intermediate degrees are
also given by the geometric definition in the introduction; in particular, if
$m=n$, then $\delta_n(P)$ is the {\it topological degree} of $P$, which is
defined as the cardinality of $P^{-1}(x)$, for a generic point
$x\in\PP^n$. Clearly, $\delta_\ell (P)>0$ if and only if $\rank\:P\geq \ell$.
It is easy to verify that $\delta_1(P) = \deg P$, and if $P$ is holomorphic
(this can happen only if $m\geq\rank\:P=n$), then $\delta_\ell (P)=(\deg
P)^\ell$, for $\ell\leq n$.  

It follows from Lemma~\ref{general formula}
applied to the maps
$P_1=\cdots=P_\ell=P\,,\,P_{\ell+1}=\cdots=P_n=\Id_{\PP^n}$ that in general,
\begin{equation}\delta_\ell (P) \leq (\deg P)^\ell\;,\quad\mbox{\rm for}\ \ 
1\leq\ell\leq
n\;.\label{delta-l}\end{equation} We shall give a more general inequality in
Lemma~\ref{k+l} below.

\begin{lemma}  Let $P:\PP^n\merto \PP^m$ be a rational
map.  Then $$\delta_\ell(P)=\deg P^{-1}(W)\leq(\deg P)^\ell$$ for generic
$W\in\glm$, with equality if and only if $\codim I_P>\ell$. In particular,
$\delta_n (P)\leq (\deg P)^n$, with equality if and only if $P$ is holomorphic.
\label{algebra}\end{lemma}

\begin{proof} We shall apply Lemma~\ref{general
formula} with
$P_1=\cdots=P_\ell=P\,,\,P_{\ell+1}=\cdots=P_n=\Id_{\PP^n}$.  By part~(ii) of
the lemma, $\delta_\ell(P)\leq(\deg P)^\ell$. For
generic $W=H_1\cap\cdots\cap H_\ell$, $\pi_2^{-1}(W)$ is of codimension $\ell$
and has no components contained in the exceptional locus of $\pi_1$, and thus
$P^{-1}(W)$ has pure dimension
$n-\ell$ and $\dim P^{-1}(W)\cap I_P<n-\ell$.  Hence for generic
hyperplanes $H_1,\ldots,H_n$, we have
$$\bigcap_{j=1}^n P_j^{-1}(H_j)\setminus I_P=P^{-1}(W)\cap
H_{\ell+1}\cap\cdots\cap H_n\;,$$ where $W=H_1\cap\cdots\cap H_\ell$, and thus
by part~(i),
$$\delta_\ell(P)=\#\left(P^{-1}(W)\cap H_{\ell+1}\cap\cdots\cap
H_n\right)=\deg P^{-1}(W)\;.$$ Furthermore, using the notation in the proof of
Lemma~\ref{general formula},
$$\{z\in\PP^n: Q_1(z)=\cdots=Q_\ell(z)=0\}=P^{-1}(W)\cup I_P\;.$$ Since $\deg
Q_1=\cdots=\deg Q_\ell=\deg P$, it follows from B\'ezout's theorem (see for
example \cite[Example 8.4.6]{Fu}) that $\delta_\ell(P)<(\deg P)^\ell$ if $\dim
I_P\geq n-\ell$. If $\dim I_P<n-\ell$, then $P^{-1}(W)\supset I_P$ (by
dimension considerations) and B\'ezout's theorem gives
equality.\qed\end{proof}

In fact, if $n=2$ in Lemma~\ref{algebra}, then $\delta_2(P)=(\deg P)^2-q$,
where $q$ is the number of points of $I_P$ counting multiplicity.  This is
illustrated by the following example.

\medskip
\begin{example}  Let $P:\PP^2 \merto \PP^2$ be
given by
$$P(z_0, z_1, z_2) = (z_1 z_2, z_0 z_2, z_0 z_1) =
\left(\frac{1}{z_0}, \frac{1}{z_1},
\frac{1}{z_2}\right)\;.$$
Then $\delta_1(P) = 2,\ \delta_2 (P) = 1$.  Note that in this 
example $I_P$ consists of the three points $(1,0,0), (0,1,0), (0,0,1)$.
\end{example}

\begin{lemma}
Suppose that $P:\PP^n \merto \PP^m$ is a 
rational
map and $L:\C^{m+1}\to\C^{M+1}$ is a linear map such that $\Image
P\not\subset\PP(L^{-1}(0))$.  Let $P_L=\wh L\circ P:\PP^n\merto\PP^M$, where
$\wh L:\PP^m\merto\PP^M$ is the map induced from $L$.  Then
$$\int_{\PP^n\setminus I}(P_L^*\omega)^k\wedge(P^*\omega)^{\ell-k}
\wedge\omega^{n-\ell}\leq\delta_\ell(P)$$ for $1\leq k\leq\ell\leq n$,
where $I=I_P\cup I_{P_L}$.\label{mixed}\end{lemma}

\begin{proof} We can assume without loss of generality that $M=m$.  We
first consider the case where $L$ is nonsingular and thus $\wh L$ is
biholomorphic.  Therefore
$$\int_{\PP^n\setminus I}(P_L^*\omega)^k\wedge(P^*\omega)^{\ell-k}
\wedge\omega^{n-\ell}=\int_{\PP^n\setminus
I}(P^*\omega')^k\wedge(P^*\omega)^{\ell-k}
\wedge\omega^{n-\ell}$$ where $\omega'=\wh L^*\omega$. Since $\omega'$ and
$\omega$ are in the same de Rham class, it follows from 
Section~\ref{s-notation} (or by the proof of Lemma~\ref{topology}) that
$$\int_{\PP^n\setminus
I}(P^*\omega')^k\wedge(P^*\omega)^{\ell-k}
\wedge\omega^{n-\ell}=\int_{\PP^n\setminus
I}(P^*\omega)^{\ell}\wedge\omega^{n-\ell}=\delta_\ell(P)\;.$$

We now suppose that $L$ is singular. Choose a sequence $\{L_\nu\}$ of
nonsingular linear operators on $\C^{m+1}$ such that $L_\nu\to L$.  We
can write
$$(P_{L_\nu}^*\omega)^k\wedge(P^*\omega)^{\ell-k}
\wedge\omega^{n-\ell}=f_\nu\omega^n\;,\quad
(P_L^*\omega)^k\wedge(P^*\omega)^{\ell-k}
\wedge\omega^{n-\ell}=f\omega^n$$ 
where $f_\nu,f$ are non-negative $\cC^\infty$ functions on $\PP^n\setminus
I$.  Then $f_\nu\to f$ pointwise on $\PP^n\setminus I$, and hence by
Fatou's Lemma,
$$\int_{\PP^n\setminus I}(P_L^*\omega)^k\wedge(P^*\omega)^{\ell-k}
\wedge\omega^{n-\ell}=\int_{\PP^n\setminus
I}f\omega^n \leq\liminf_{\nu\to\infty}\int_{\PP^n\setminus I}
f_\nu\omega^n =\delta_\ell(P)\;.$$
\qed\end{proof}

\begin{lemma} Let $P:\PP^n\merto \PP^m\;,\;Q:\PP^m\merto \PP^r$ be rational
maps.  Then $$\delta_\ell(Q\circ P)\leq\delta_\ell(P)\delta_\ell(Q)\;.$$
\label{qp}\end{lemma}

\begin{proof} Let $\eta=Q^*\omega^\ell\in\cD'{}^{\ell,\ell}(\PP^n)$.  We
smooth $\eta$ by an approximate identity $\{\psi_\epsilon\}_{\epsilon>0}$
with respect to a Haar measure $h$ on $GL(n+1,\C)$ to obtain
$$\eta_\epsilon\eqd\int\limits_{GL(n+1,\C)}(g^*\eta) \psi_\epsilon(g)dh(g)
\in\cE^{\ell,\ell}(\PP^n)\;.$$  Then $\eta_\epsilon\to\eta$ pointwise as
$\epsilon\to 0$, $\eta_\epsilon\geq 0$, and we have the identity in de Rham
cohomology,
$$[\eta_\epsilon]=[\eta]=\delta_\ell(Q)[\omega^\ell]\in
H^{2\ell}(\PP^n,\Z)\;.$$  Using the commutative diagram (\ref{comm-diagram}),
we then have $$\int_{\PP^n}P^*\eta_\epsilon\wedge\omega^{n-\ell}=\int_{\wt
G}\pi^*_2\eta_\epsilon\wedge\pi_1^*\omega^{n-\ell}=\int_{\wt G}
\delta_\ell(Q)\pi_2^*\omega^\ell\wedge\pi_1^*\omega^{n-\ell}=\delta_\ell(Q)
\delta_\ell(P)\;.$$ Therefore, by Fatou's lemma,
$$\delta_\ell(Q\circ P)=\int_{\PP^n\setminus I_{Q\circ
P}}P^*\eta\wedge\omega^{n-\ell}\leq\liminf_{\epsilon\to 0}
\int_{\PP^n}P^*\eta_\epsilon\wedge\omega^{n-\ell}=\delta_\ell(Q)
\delta_\ell(P)\;.$$
\qed\end{proof} 

\begin{lemma} Let $P:\PP^n\merto \PP^m$ be a rational
map.  Then $$\delta_{k+\ell}(P)\leq\delta_k(P)\delta_\ell(P)\;.$$
\label{k+l}\end{lemma}

\begin{proof}  Let $\eta=P^*\omega^\ell\in\cD'{}^{\ell,\ell}(\PP^m)$ and
consider the smooth forms $\eta_\epsilon$ as in the above proof.  As
before,
$$[\eta_\epsilon]=\delta_\ell(P)[\omega^\ell]\in
H^{2\ell}(\PP^n,\Z)\;,$$ and
\begin{eqnarray*}\int_{\PP^n}P^*\omega^k\wedge\eta_\epsilon\wedge
\omega^{n-k-\ell}&=&\int_{\wt
G}\pi^*_2\omega^k\wedge\pi_1^*(\eta_\epsilon\wedge
\omega^{n-k-\ell})\\&=&\delta_\ell(P)\int_{\wt G}
\pi_2^*\omega^k\wedge\pi_1^*\omega^{n-k}=\delta_\ell(P)
\delta_k(P)\;.\end{eqnarray*}  The conclusion follows as above by letting
$\epsilon\to 0$ and applying Fatou's lemma.\qed\end{proof}

%%%%%%%%%%%%
\section {The proximity function \label{s-proximity}}

Let $P:\PP^n \merto \PP^m$ be a rational map, and let $1\leq\ell \leq \nm$.
Recall that if $\codim P^{-1}(W) =\ell$, then equations (\ref{lambda1}),
(\ref{lambda2})
define a current $P^*\Lambda_W\in\cD'{}^{\ell-1,\ell-1}(\PP^n)$ with
locally integrable coefficients.  Hence we can define the {\it proximity
function} $m^\ell_P:\glm\to [0,+\infty]$ by
$$\hskip-1em
m^\ell_P(W)=\left\{\begin{array}{ll}(P^*\Lambda_W,\omega
^{n-\ell+1})=\int_{\PP^n\setminus(I_P\cup I_{P,W})}P^*\Lambda_W\wedge \omega
^{n-\ell+1}\quad &\mbox{\rm if}\ \codim P^{-1}(W) =\ell\\ +\infty &
\mbox{\rm if}\ \codim P^{-1}(W) < \ell\end{array} \right.\hskip-1em
$$ for $W\in\glm$,
where $I_{P,W}$ is the indeterminacy locus of the map $z\mapsto P(z)\wedge
W\in \PP(\bigwedge^{m+2-\ell}\C^{m+1})$. We give $\glm$ the K\"ahler metric
$\omega$ induced from the natural embedding $\glm\subset
\PP(\bigwedge^{m+1-\ell} \C^{m+1})$.

The following key estimate is used in our proof of Theorem~\ref{th1}.

\begin{lemma}
If $P:\PP^n \merto \PP^m$ is a
rational map, then $m^\ell_P\in {\cal L} ^1
(\glm)$ and 
$$dd^cm^\ell_P\leq \ell \delta_{\ell-1}(P) \omega \;,$$
for $1\leq\ell\leq \nm$.\label{ddcmp}
\end{lemma}

\begin{proof} (The techniques used for this proof are somewhat similar to 
those of \cite{Skoda-1,Skoda-2}.) 
Assume first that $\ell\geq 2$. (The
estimate of Lemma~\ref{ddcmp} is straightforward for the hyperplane case
$\ell=1$; we give the argument for this case at the end of this proof. 
See also [RS] for a complete treatment of pre-images of hyperplanes.) 
Write $E =
\bigwedge^{m+1-\ell}\C^{m+1}$, $\tilde E=\bigwedge^{m+2-\ell}\C^{m+1}$
and let $\lambda _0:E\rightarrow \C,\ \lambda_1:\tilde E\rightarrow \C $
be linear functions of unit norm.  For $\zeta \in \C^{m+1}$, $W\in E$, we
define the {\it augmented exterior product}
$$\zeta \tilde{\wedge} W = (\zeta_0 W_0, \zeta\wedge W)\in \C\oplus
\tilde E$$ where $W_0 = \lambda _0 (W)$.  Let
$$X = \PP^n_z \times \glm _W\times \PP^m_\zeta
\times \PP(\tilde E)_\theta \times \PP(\C\oplus
\tilde E)_\eta\;.$$ (The subscripts $z$, $W$, $\zeta$, $\theta$, $\eta$
serve to identify the variables used in this discussion.)
\end{proof}

We further write $\lambda _1(\theta)=\theta_1$ for
$\theta\in\tilde E$.  If $\eta = (c,\theta)\in
\C\oplus \tilde E$, we write $\eta_0=c$,
$\eta_1=\theta_1=\lambda _1(\theta)$.  By making a
linear change of coordinates in $\PP^m_\zeta$ we
can assume without loss of generality that
Image$(P) \not\subset \{\zeta\in
\PP^m:\zeta_0=0\}$. We consider the current
$$\Omega =\log
\frac{|\zeta|^2|W|^2|\theta_1|^2|\eta_0|^2}
{|\zeta_0|^2|W_0|^2|\theta|^2|\eta_1|^2}
\omega
_z^{n-\ell+1}\wedge\sum^{\ell-1}_{j=0}\omega^
j_\theta\wedge \omega ^{\ell-1-j}_\zeta \in
\cD'{}^{n-1,n-1}(X)\;.$$
Let $Q:\PP^n\times \glm \merto X$ be the
meromorphic (rational) map given by
\begin{equation} Q(z, W) = (z, W, P(z), P(z)\wedge W, P(z)\tilde
\wedge W)\;.\label{q}\end{equation}
(The redundancy in the map $Q$ is needed
to facilitate our proof.)
We can assume that $\lambda_1 $ was chosen so that
$\lambda_ 1 (P(z)\wedge W)\not\equiv 0$.  We shall
show (see Lemma~\ref{m}) that
$$m_P^\ell = \pi_{2*} Q^*\Omega\;,$$
where $\pi_2:\PP^n\times \glm \rightarrow \glm $
is the projection.  Of course, since $Q$ is not
regular, we must define the pull-back $Q^*\Omega$, which we do as
follows.  Let
$$Y\stackrel{\rho}{\rightarrow } \operatorname{Image}(Q)\subset X$$
be a desingularization of the image of $Q$.  (Note
that Image$(Q)$ is an algebraic subvariety of $X$
and can be identified with the graph of $Q$.)  Let
$\pi_1:X\rightarrow \PP^n\times \glm $ be the
projection, and write $\rho_1 = \pi_1\circ \rho,
\rho_2=\pi_2\circ \rho_1$ so that we have the
commutative diagram:
$$\begin{array}{llll}
\PP^n\!\times\!\glm \;\stackrel{\rho_1}{\longleftarrow}\;Y\\[.1cm] 
\hspace*{.6cm}Q\;\stackrel{|}{\downarrow}\;\rho\swarrow\;\;
\;\;\;\;\downarrow\;\rho_ 1\;\;\; 
\;\;\searrow\rho_2\\ [.2cm]
\hspace*{1.1cm}X\stackrel{\pi_1}{\longrightarrow}\;\PP^n\!\times\!
\glm \stackrel{\pi_2}{\rightarrow} \glm 
\end{array}
$$
We then define $Q^*\Omega = \rho_{1*}\rho^*\Omega$, where $\rho^*\Omega$
is given as in Section~\ref{s-notation}.  Let $I_Q\subset 
\PP^n\times
\glm $ denote the indeterminacy locus of $Q$; $I_Q$ is an algebraic
subvariety of codimension $\geq 2$.  Write $U = \PP^n\times
\glm \setminus I_Q$, and let
$Q_0=Q|_U:U\rightarrow X$.  Since $\rho_1$ maps $\rho_1^{-1}(U)$
biholomorphically onto $U$, we have $Q_0=\rho\circ\rho_1^{-1}|_U$. 
Hence
$$(Q^*\Omega)|_U=(\rho_1^{-1}|_U)^*\rho^*\Omega = Q^*_0\Omega\;,$$
so our definition of $Q^*$ agrees with the usual one at regular
points.  Since $Q^*\Omega$ has coefficients in ${\cal L}^1_{\loc}$,
it is the extension to $\PP^n\times \glm $ of $Q^*_0\Omega$ with zero mass on
$I_Q$.

We first note the following:

\begin{lemma}
$m_P^\ell=\rho_{2*}\rho^*\Omega=\pi_{2*}Q^*\Omega \in {\cal L}^1(\glm )$.
\label{m}\end{lemma}

\begin{proof}  By the definition of $Q^*\Omega$ we have 
\begin{equation}\rho_{2*} \rho^*\Omega = \pi_{2*}\rho_{1*}\rho^*\Omega = \pi_{2*}
Q^*\Omega\;.\label{mp1}\end{equation} Since $Q^*\Omega$ has coefficients in ${\cal
L}^1_{\loc}$ and its degree, $2n$, is the fibre dimension of $\pi_2$,
$\pi_{2*}Q^*\Omega\in {\cal L}^1(\glm )$ and \begin{equation}(\pi_{2*}Q^*\Omega)(W) =
\int_{\PP^{n}_{z}\times \{W\}\setminus I_Q}Q^*\Omega<+\infty\;,\quad\mbox{\rm
for a.a. } \ W\in \glm \;.\label{mp2}\end{equation} 
Since $\eta_1\circ Q =
\theta_1\circ Q$ and $\eta_0\circ Q = P_0(z)W_0=(\zeta_0W_0)\circ Q$, we have
$$Q^*\Omega|_U = Q_0^*\Omega=\log\frac{|P(z)|^2|W|^2}{|P(z)\wedge
W|^2}\omega ^{n-\ell+1}_z\wedge \sum^{\ell-1}_{j=1}
(dd^c\log|P(z)\wedge W|^2)^j \wedge (P^*\omega _\zeta)^{\ell-1-j}\;.$$
Thus for a.a. $W\in\glm $, $I_Q\not\supset \PP^n\times \{W\}$ and
\begin{equation} m_P^\ell(W)=\int_{\PP^n\setminus(I_P\cup I_{P,W})}P^*\Lambda_W\wedge
\omega ^{n-\ell+1}=\int_{\PP^n_z\times \{W\}\setminus I_Q}Q^*\Omega\leq
+\infty\;. 
\label{mp3}\end{equation} The desired identity follows from (\ref{mp1}), (\ref{mp2}),
and(\ref{mp3}). \hfill$\Box$
\end{proof}

\begin{rem} We could use (\ref{lambda1}), (\ref{lambda2}) to 
define $m_P^\ell(W)$ for all
$W\in \PP(\bigwedge^{m+1-\ell} \C^{m+1})$; then Lemma~\ref{ddcmp} remains valid
on $\PP(\bigwedge^{m+1-\ell} \C^{m+1})$.
\end{rem}

We are now ready to compute $dd^cm_P^\ell$.  We let $H^0_\zeta$, $H^0_W$,
$H^1_\theta$, $H^0_\eta$, $H^1_\eta$ denote the hyperplanes in $X$ given by
the divisors of $\zeta_0,W_0,\theta_1, \eta_0,\eta_1$ respectively, and we let
$$D = H^0_\zeta + H^0_W-H^1_\theta+H^1_\eta - H^0_\eta = {\operatorname{
Div}}\left(\frac{\zeta_0W_0\eta_1}{\theta_1\eta_0}\right)\;.$$
We have
\begin{equation}\begin{array}{ll}
dd^c\Omega &=(\omega_ \zeta +\omega _W - \omega _\theta - D)\wedge
\omega _z^{n-\ell+1} \wedge \displaystyle{\sum^{\ell-1}_{j=0}} \omega
^j_\theta
\wedge \omega _\zeta^{\ell-1-j}\\[.2cm]
&=\left[\omega ^\ell_\zeta - \omega^\ell_\theta + (\omega_ W -
D)\wedge \displaystyle{\sum^{\ell-1}_{j=0}}\omega ^j_\theta\wedge
\omega ^{\ell-1-j}_\zeta \right]\wedge \omega
_z^{n-\ell+1}\;.\end{array}\label{ddcomega}\end{equation}

\begin{lemma}
The divisor $\rho^*D$ is effective (i.e., is
locally the divisor of a holomorphic function).
\label{effective}\end{lemma}

\begin{proof}  (Our choices of $\zeta_0$ and $\theta_1$ guarantee
that $\rho(Y)\not\subset \Supp\:D$ so that $\rho^*D$ is defined.) 
Let $y_0\in Y$ be arbitrary, and let 
$$X'=\{(z, W, \zeta, \theta, \eta) \in X: z_i\neq 0, \lambda
(W)\neq 0, \zeta _j\neq 0, \lambda '(\theta)\neq 0\}$$
where $0\leq i\leq n$, $0\leq j\leq m$ and $\lambda :E\rightarrow
\C, \lambda ':\tilde E \rightarrow \C$ are linear functions chosen
so that $\rho (y_0)\in X'$.  Let
$$g = \frac{\tilde \zeta_0\tilde W_0}{\tilde
\theta_1}\frac{\eta_1}{\eta_0}\in {\operatorname{Mer}}(X')$$
where $\tilde \zeta_0 = \zeta_0/\zeta_j$, $\tilde W_0 = W_0
/\lambda(W) = (\lambda _0/\lambda)(W)$, $\tilde \theta_1 =
\theta_1/\lambda '(\theta) = (\lambda _1/\lambda ')(\theta)$. 
Then Div$(g) = D|_{X'}$.  We must show  that $g\circ \rho$ is
holomorphic at $y_0$.  Now
$$\begin{array}{ll}
\displaystyle{g\circ Q = \frac{P_0(z)}{P_j(z)}\frac{W_0}{\lambda
(W)}\frac{\lambda '(P(z)\wedge W)}{\lambda _1(P(z)\wedge
W)}\frac{\lambda_ 1(P(z)\wedge W)}{P_0(z)W_0} }\\[.5cm]
\hspace*{2cm}=\displaystyle{ \frac{\lambda'
(P(z)\wedge W)}{P_j(z) \lambda (W)} = \frac{\lambda '(P(\tilde
z)\wedge \tilde W)}{P_j(\tilde z)} }\end{array}$$
where $\tilde z = z^{-1}_i z\in\C^{n+1}$, $\tilde W = \lambda (W)
^{-1}W\in E$.  Write, for $y \in Y$, 
$$\rho_1(y) = (z(y), W(y))\;,$$
$\tilde z(y) = z^{-1}_i(y)z(y), \tilde W(y) = \lambda (W(y))^{-1}
W(y)$.  Since $\rho$ is holomorphic and $\rho(y_0) \not\in$
Div$(\zeta_j)$, there is a neighborhood $Y_0$ of $y_0$ so that
$P_\mu (\tilde z(y)) = \varphi(y) f_\mu (y)$ for $y\in Y_0$ and $0\leq \mu
\leq m$, where $\varphi, f_0,\ldots, f_m \in {\cal O}(Y_0), 
\varphi \not\equiv 0$, and $f_j(y_0)\neq 0$.  Thus
$$g\circ\rho = g\circ Q\circ \rho_1 = \frac{\lambda '(\varphi
F\wedge\tilde W)}{\varphi f_j} = \frac{\lambda '(F\wedge\tilde
W)}{f_j}$$
on $Y_0$, where $F = (f_0,\ldots f_m)$.  Since $f_j(y_0)\neq 0$,
$g\circ \rho$ is holomorphic at $y_0$. \hfill $\Box$
\end{proof}

\medskip We now complete the proof of Lemma~\ref{ddcmp}:  By
Lemma~\ref{effective} and (\ref{ddcomega}), we have
$$\rho^*dd^c\Omega \leq \rho^*\left[\left(\omega ^\ell_\zeta +\omega
_W
\wedge
\sum^{\ell-1}_0 \omega ^j_\theta\wedge
\omega_\zeta^{\ell-1-j}\right)\wedge
\omega ^{n-\ell+1}_z\right]\;.$$
By Section~\ref{s-notation}, $\rho^*dd^c\Omega = dd^c\rho^*\Omega$.  Hence by
Lemma~\ref{m},
$$\begin{array}{ll}
dd^cm_P^\ell&=dd^c(\rho_{2*}\rho^*\Omega) = \rho_{2*} \rho^*
dd^c\Omega\\[.2cm]
&\leq \rho_{2*} \rho^*\left[(\omega^\ell_\zeta +\omega _W\wedge
\sum^{\ell-1}_0 \omega ^j_\theta\wedge \omega_
\zeta^{\ell-1-j})\wedge \omega _z^{n-\ell+1}\right]\;.
\end{array}$$
Since $\rho_{2*}\rho^* = \pi_{2*} Q^*$ and
$$Q^*_0(\omega ^\ell_\zeta \wedge \omega _z^{n-\ell+1}) = (dd^c\log
|P(z)|^2)^\ell \wedge \omega _z^{n-\ell+1}=0\;,$$
we obtain
$$\begin{array}{ll} dd^cm_P^\ell&\leq \displaystyle{\sum^{\ell-1}_{j=0}
\pi_{2*} Q^*\left(\omega _W\wedge \omega ^j_\theta\wedge \omega
_\zeta^{\ell-1-j}\wedge \omega _z^{n-\ell+1}\right)}\\[.5cm]
&=\displaystyle{\sum^{\ell-1}_{j=0} \left(\int_{\PP^n_z\times \{W\}} (dd^c\log
|P(z)\wedge W|^2)^j\wedge(P^*\omega _\zeta)^{\ell-1-j}\wedge \omega
_z^{n-\ell+1}\right)\omega _W}\;.\end{array}$$ 
By Lemma~\ref{mixed} with
$L:\C^{m+1}\to \wt E$ given by $L(\zeta)=\zeta\wedge W$,
$$\int_{\PP^n_z\times \{W\}}(dd^c\log |P(z)\wedge W|^2)^j\wedge
(P^* \omega_{\zeta})^{\ell-1-j}\wedge \omega_z^{n-\ell+1} \leq
\delta_{\ell-1}(P)\;,$$
and the desired inequality follows.

We now modify (and simplify) the above argument for the case $\ell=1$:
Identify
$$G(1,m)=\PP(\textstyle{\bigwedge^m}\C^{m+1})=\PP(\C^{m+1*})=\PP^{m*}\;.$$
Let
$X=\PP^n_z\times\PP^{m*}_W\times\PP^m_\zeta$ and consider the current
$$\Omega=\log\frac{|\zeta|^2
|W|^2}{|(W,\zeta)|^2}\omega^n_z\in\cD'{}^{n,n}(X)\;.$$ Then
$$\ddc\Omega=(\omega_\zeta+\omega_W-D)\wedge\omega^n_z\;,$$
where $D={\operatorname{Div}}(W,\zeta)$.  Then $m_P^1=\pi_{2*}Q^*\Omega$,
where $$Q=\left(z,W,P(z)\right):\PP^n_z\times\PP^{m*}_W\to X\;.$$
We conclude as
before that $\ddc m_P^1=\pi_{2*}Q^*\ddc\Omega\leq\omega_W$.
\hfill $\Box$

%%%%%%%%%%%%%%
\section {Description of the exceptional set \label{s-exceptional}}

In this section we prove Theorem~\ref{th1}, giving a description of the
exceptional set $\cE$ in terms of the proximity function as follows.  Let $\cP
= \{P_k\}$ be a sequence of rational mappings from $\PP^n$ to $\PP^m$ as in
Theorem~\ref{th1} and fix a sequence $\cA=\{a_k\}$ of positive numbers.  We
define the {\it exceptional set} \begin{equation}\cE_\ell(\cP,\cA)=\left\{W\in
\G(\ell,m):\limsup_{k\rightarrow
\infty}\frac{m_k(W)}{a_k}>0\right\}\;,\label{exceptional}\end{equation} where
$m_k=m_{P_k}^\ell$. (Thus by definition, $\cE_\ell(\cP,\cA)$ contains those
planes $W$ such that infinitely many of the pre-images $P_k^{-1}(W)$ have
codimension less than $\ell$.)  The following two propositions yield
Theorem~\ref{th1} with $\cE=\cE_\ell(\cP,\cA)$.

\begin{proposition}  Let $\cP=\{P_k\}$ be a sequence of rational mappings
from $\PP^n$ to $\PP^m$ and let $\cA=\{a_k\}$ be a sequence of positive 
numbers. Let $1\leq\ell \leq \nm$. 
Then for all $W\in \G(\ell, m)\setminus
\cE_\ell (\cP,\cA)$,
$$\frac{1}{a_k}\left(P^*_k
[W]-P^*_k\omega^\ell\right)\rightarrow 0$$ as $k\rightarrow \infty$.
\label{pr-i}\end{proposition}

\begin{proposition} Let $\cP=\{P_k\},\;\cA=\{a_k\}$ be as in 
Proposition~\ref{pr-i}.  If $$\sum_{k=1}^\infty\frac{\delta_{\ell-1}
(P_k)}{a_k}<+\infty\;,$$ then $\cE_\ell(\cP,\cA)$ is pluripolar in
$\G(\ell,m)$.
\label{pr-ii}\end{proposition}

Before proving Propositions \ref{pr-i} and \ref{pr-ii}, we note a corollary to
Theorem~\ref{th1} on the equidistribution of pre-images for subsequences of a
given sequence of rational mappings.  This corollary uses the following
terminology:  For a current $T\in\cD'{}^{p,p}(\PP^n)$ of
order 0, we let $\|T\|$ denote the total variation measure of $\PP^n$, which is
the regular measure on $\PP^n$ given by
$$\|T\|(U)=\sup\left\{|(T,\phi)|:\phi\in\cD^{n-p,n-p}(U), \|\phi\|\leq
1\right\}$$ for $U$ open in $\PP^n$.  Here, $\|\phi\|$ denotes the comass norm
of a compactly supported form $\phi$ (see Federer \cite[1.8.1, 4.1.7]{Fe}).
The quantity $\|T\|(\PP^n)$ is called the mass of T. If $T$ is positive, then
it follows from Wirtinger's inequality (see for example, \cite[1.8.2]{Fe})
that $\|T\|(\PP^n)=(T,\omega^{n-p})$. In particular, $\|\omega^p\|(\PP^n)=1$.

\begin{cor} Let $\cP=\{P_k\}$ be as in Theorem~\ref{th1}, let
$1\leq\ell \leq \nm$, and suppose that $$\frac{\delta_{\ell-1}
(P_k)}{\delta_\ell (P_k)}\rightarrow 0\;.$$  Let $\cM \subset
\cD'{}^{\ell, \ell}(\PP^n)$ be the set of weak cluster points of
$\{\frac{1}{\delta_\ell (P_k)} P_k^*\omega ^\ell\}$.  Then $\cM\neq
\emptyset$, and for every current $\eta \in \cM$,
\begin{itemize}
\item[\rm i)] $\eta$ is a positive current
\item[\rm ii)]$\eta$ has mass $1$
\item[\rm iii)] there is a subsequence $\{P'_k\}$ of
$\{P_k\}$ and a pluripolar set $\cE$ such that 
$$\frac{1}{\delta_\ell (P'_k)}P'{}^*_k [W]\rightarrow \eta$$ for all $W\in
\G(\ell, m)\setminus \cE$.
\end{itemize}
\label{cor2}\end{cor}

\noindent{\it Proof (assuming Theorem~\ref{th1}):\/} Since
$$\|\frac{1}{\delta_\ell (P_k)} P^*_k\omega
^\ell\|(\PP^n)=\frac{1}{\delta_\ell(P_k)}\int P^*_k\omega^ \ell\wedge \omega
^{n-\ell}=1,$$ it follows that $\cM\neq \emptyset$.  Let
$$\eta = \lim_{k\to\infty}\frac{1}{\delta_\ell(P'_k)} P^{'*}_k \omega
^\ell\in\cM$$ for some subsequence $\{P'_k\}$ of $\cP$.  Then (i) is obvious,
and (ii) follows from
$$\|\eta\|(\PP^n)=(\eta, \omega ^{\ell -1}
)=\lim_{k\to\infty}\frac{1}{\delta_\ell(P'_k)}(P^{'*}_k\omega ^\ell, \omega
^{n-\ell})=1\;.$$ Finally, choose a subsequence $\{P''_k\}$ of $\{P'_k\}$ such
that $\displaystyle\sum\frac{\delta_{\ell-1}(P''_k)}{\delta_\ell(P''_k)}<+
\infty $ , and apply Theorem~\ref{th1} with $a_k=\delta_\ell(P''_k)$ to obtain
(iii).  \qed

\smallskip
\noindent{\it Proof of Proposition \ref{pr-i}:\/}  Let 
$\phi\in\cD^{n-\ell,n-\ell}(\PP^n)$ be an arbitrary real form and choose a
constant $c_\phi$ such that $-c_\phi\omega^{n-\ell+1}\leq dd^c\phi\leq
c_\phi\omega^{n-\ell+1}$. Then for all $W\in\glm$ such that $\codim
P_k^{-1}(W)=\ell$, we have by (\ref{usual-Poincare-Lelong}),
$$|(P_k^*[W] -P_k^* \omega^{\ell}, \phi)| = |(P_k^* \Lambda_W, \ddc 
\phi)| \leq c_\phi m_k(W)\;.$$
The conclusion follows 
from the definition of $\cE_\ell(\cP,\cA)$.\qed

\smallskip
\noindent{\it Proof of Proposition \ref{pr-ii}:\/}
Let $P:\PP^n\merto\PP^m$ be a rational map.  We write
$$S_P=\left\{W\in \G(\ell,m):\codim P^{-1}(W)<\ell\right\}\;;$$ $S_P$ is an
algebraic subvariety (which is usually empty) of $\G(\ell,m)$ and thus is
pluripolar.  By (\ref{mp3}) we can write
$$P^*\Lambda_W\wedge\omega_z^{n-\ell+1} = f(z,W)\omega_z^n\;,$$ where
\begin{itemize}\item[i)]

$0\leq f(z,W)\leq+\infty$,
\item[ii)]
$f\in\cL^1(\PP^n\times \G(\ell,m))$,
\item[iii)] $f$ is continuous on the
set $U\eqd \PP^n\times \G(\ell,m)\setminus I_Q\;.$ 
\end{itemize}  Here $Q$ is the map given by equation~(\ref{q}) in
Section~\ref{s-proximity}.  By the definition of the proximity function,
$$m_P^\ell(W)=\int_{\PP^n_z}f(z,W)\omega^n_z\;,\quad\mbox{\rm for}\ \ \
W\in\glm\setminus S_P\;.$$  Since
$$I_Q\cap(\PP^n\times\{W\})=\left[I_P\cup P^{-1}(W)\right]\times \{W\}\;,$$
$m_P^\ell(W)\in[0,+\infty]$ is well defined by the above integral.

Let $(-m^\ell_P)^*$ be the upper-regularization of $(-m^\ell_P)$ given by
$$(-m^\ell_P)^*(W)=\limsup_{\xi\to W}[-m^\ell_P(\xi)]\;,$$ for $W\in\glm$, and
write $$m^\ell_{P*}(W) =-(-m^\ell_P)^*(W)=\liminf_{\xi\to W}m^\ell_P(\xi)\;.$$
We claim that 
\begin{equation} m^\ell_P(W)= m^\ell_{P*}(W) \ \ \mbox{\rm for all} \
W\in\glm\setminus S_P\;.\label{lowerreg}\end{equation} 
To demonstrate (\ref{lowerreg}), let $W_0\in\glm\setminus S_P$, and 
let $\{W_\nu\}$ be a sequence of points converging to
$W_0$ such that $m^\ell_P(W_\nu)\to m^\ell_{P*}(W_0)$. Since $f(z,W_\nu)\to
f(z,W_0)$ for almost all
$z\in\PP^n$ (precisely, for $z\not\in I_P\cup P^{-1}(W_0)$), we have by
Fatou's Lemma,
$$m^\ell_P(W_0)=\int_{\PP^n_z}f(z,W_0)\omega^n_z\leq \liminf_{\nu\to\infty}
\int_{\PP^n_z}f(z,W_\nu)\omega^n_z=\liminf_{\nu\to\infty}m^\ell_P(W_\nu)
=m^\ell_{P*}(W_0)\;.$$ Since by definition, $m^\ell_P(W_0)\geq
m^\ell_{P*}(W_0)$, we obtain  (\ref{lowerreg}).

A function $u$ on a complex manifold $Y$ (with values in $\R\cup
-\{\infty\}$) is said to be {\it quasi-plurisubharmonic\/}
({\it quasi-psh\/} for short) if $u$ is locally equal to the sum of a
$\cC^\infty$ 
function and of a plurisubharmonic (psh) function, or equivalently, if
$\ddc u$ is bounded below by a continuous real $(1,1)$-form.
By Lemma~\ref{ddcmp},
$-m^\ell_{P*}$ is quasi-psh.

We now let $\cP=\{P_k\}$, $\cA=\{a_k\}$
be as in Proposition~\ref{pr-ii} and suppose that 
$\sum\frac{\delta_{\ell-1}(P_k)}{a_k}<+\infty$. It suffices to show that
$\cE_\ell(\cP,\cA)$ is pluripolar in an arbitrary affine  open set
$$G'\eqd\{W\in\glm:\lambda(W)\neq 0\}$$ given by a linear function
$\lambda:E\to\C$. . We write
$m_k=m_{P_k}^\ell$ and we  
consider the functions $u_k:G'\to\R\cup\{-\infty\}$
given by 
$$u_k = v -\frac{m_{k*}}{\ell \delta_{\ell-1}(P_k)}$$ where \begin{equation} v(W)=\log 
\frac{|W|^2}{|\lambda(W)|^2}\;.\label{v}\end{equation} We assume that $\|\lambda\|=1$, so
$v(W)\geq 0$. By Lemma~\ref{ddcmp},
$u_k$ is plurisubharmonic on $G'$. Next, construct the series 
\begin{equation} u=\sum \limits_{k=1}^\infty 
\frac{\delta_{\ell-1}(P_k)}{a_k}u_k\;.\label{u}\end{equation} 
To see that $u$ is
plurisubharmonic, represent it as a limit of
the sequence 
$$ \tau_k=S\cdot v - \frac{1}{\ell}\sum\limits_{j=1}^k
\frac{m_{j*}}{a_j}\;,$$ where
$$S=\sum\limits_{k=1}^\infty 
\frac{\delta_{\ell-1}(P_k)}{a_k}\;.$$
Since
$m_{k*}\ge 0$, $\{\tau_k\}$  is a decreasing sequence of plurisubharmonic
functions on $G'$, so the limit is  either plurisubharmonic or identically
$-\infty$. To see that the latter  case is impossible, we average
$u(W)$ over all
$W$ with respect to Haar probability measure $\sigma$ on $\glm.$
It is well known (e.g., see  \cite[Ch.~2, \S~4, Th.~7]{Shabat}) that
\begin{equation}\int\limits_{\glm} \Lambda_W(z) d\sigma(W) =
c_{\ell,m} \omega^{\ell-1} (z)\label{average}\end{equation}
for some constant $c_{\ell,m}$.
Therefore, by Fubini's theorem,
\begin{eqnarray*}\int\limits_{\glm} m_{P_k}(W) d\sigma (W) &=&
\int\limits_{\glm}\left(\int\limits_{\PP^n} P_k^* \Lambda_W
\wedge \omega^{n-l+1}\right) d\sigma (W)\\
&=&c_{\ell,m} \int\limits_{\PP^n}P_k^* \omega^{l-1} \wedge
\omega^{n-\ell+1}= c_{\ell,m}\delta_{\ell-1}(P_k)\;.\end{eqnarray*}

Hence 
\begin{eqnarray*}\int u(W)d\sigma(W) &\geq & - 
\frac{1}{\ell}\sum_{k=1}^\infty \frac{1}{a_k}\int m_k(W)d\sigma(W)\\
& = & - \frac{1}{\ell}\sum_{k=1}^\infty \frac{c_{\ell,m}
\delta_{\ell-1}(P_k)}{a_k}\\
& = & -\frac{c_{\ell,m}}{\ell}S>-\infty\;.
\end{eqnarray*} 
Thus $u$ is plurisubharmonic.

Finally, if $W\in \cE_\ell(\cP,\cA)\setminus\bigcup_k S_{P_k}$, then by
definition, for an infinite number  of indices $k$ we have
$$\epsilon a_k<m_k(W)= m_{k*}(W)$$ for some 
$\epsilon>0$, and therefore
$$\frac{\delta_{\ell-1}(P_k)}{a_k}u_k(W)\leq
\frac{\delta_{\ell-1}(P_k)}{a_k}v(W)-\frac{\epsilon}{\ell}\;,$$ so
$u(W)=-\infty.$ Since the $S_{P_k}$ are pluripolar, it follows that the
set
$\cE_\ell(\cP,\cA)$ is pluripolar. \qed

\bigskip

%%%%%%%%%%%%%%%%%
\section {Limit measures for iterates of rational maps \label{s-proof2}}

In this section we prove Theorem~\ref{th2}.  We continue to use the notation
from Section~\ref{s-exceptional}.

Assume that $P:\PP^n\merto \PP^n$ is as in Theorem~\ref{th2} and 
write $\lambda=\delta_n(P),\ \delta=\delta_{n-1}(P).$ Let
$$h=\frac{1}{\lambda} P^* \omega^n \in \cD'{}^{n,n}(\PP^n).$$  By the definition
of the topological degree $\delta_n(P)$, $\int_{\PP^n}h=1$.

\noindent {\it Claim:\/} If $f$ is a quasi-plurisubharmonic function on
$\PP^n,$ then \begin{equation}\int_{\PP^n}fh>-\infty\;.\label{claim}\end{equation}

\smallskip\noindent (Note that $h\geq 0;$ since $f$ is bounded above, the
claim is equivalent to saying that $fh$ is $\cL^1.$) To verify the claim, we
again consider the commutative diagram (\ref{comm-diagram}). We then have
$$\int fh = \frac{1}{\lambda}\int_{\wt{G}} (f\circ \pi_1) 
\pi_2^*\omega^n > -\infty$$
since $f\circ \pi_1$ is quasi-plurisubharmonic on $\wt{G}$ and hence is 
in $\cL^1$, verifying (\ref{claim}).

Choose $a>0$ such that $\delta<a<\lambda.$ By Lemma~\ref{qp}, 
$\delta_{n-1}(P^k)\leq \delta^k$ and thus
$$\sum \frac{\delta_{n-1}(P^k)}{a^k} = S < +\infty.$$

Applying Theorem~\ref{th1} with $\ell=n,\ \cP=\{P_k\},$ where $P_k$ is the 
$k$-th iterate of $P,$ and $\cA=\{a^k\},$ we conclude that
\begin{equation}\frac{m_k(W)}{a^k}\to 0\label{e51}\end{equation}
for points $W$ of $\PP^n$ outside a pluripolar set.

Since pluripolar sets have Lebesgue measure zero, (\ref{e51}) is valid 
for a.a. $W\in \PP^n$. Let
$$U \eqd \sum\limits_{k=1}^\infty \frac{m_k}{a^k} = Sv-u,$$ where $v,u$ are
given by (\ref{v}),(\ref{u}).  By the proof of Proposition~\ref{pr-ii}, 
$u$ is psh
and hence $-U$ is quasi-plurisubharmonic. Therefore, by (\ref{claim}),
\begin{equation}
\int \limits_{\PP^n} Uh < + \infty \label{e52}
\end{equation}

Let $\phi\in \cD^0 (\PP^n)$ be arbitrary. Since 
$\displaystyle\int_{\PP^n} h=1,$ 
\begin{eqnarray*}
\frac{1}{a^k} (P_k^* h - P_k^* \omega^n, \phi)
&=&\frac{1}{a^k} \int (P_k^* [W] - P_k^* \omega^n, \phi) h(W)\\
&=&-\frac{1}{a^k} \int (P_k^* \Lambda_W, \ddc \phi) h(W)\;.
\end{eqnarray*}

As in the proof of Proposition~\ref{pr-i},
$|(P_k^* \Lambda_W, \ddc\phi)| \leq c_\phi m_k(W)$, and we conclude that
\begin{equation}|\frac{1}{a^k} (P_k^*\Lambda_W, dd^c\phi)| \leq c_\phi
U(W)\;.\label{e53}\end{equation}  

By (\ref{e51}), (\ref{e52}) and (\ref{e53}), we can let $k\rightarrow +\infty$
and apply Lebesgue's dominated convergence theorem to conclude that
\begin{equation}\frac{1}{a^k} (P_k^* h - P_k^* \omega^n, \phi) \to 0
\label{e54}
\end{equation}
as $k\rightarrow +\infty.$

We note that 
$$\frac{1}{\lambda^{k+1}} P_{k+1}^* \omega^n = \frac{1}{\lambda^{k}} 
P_{k}^* h,$$
as an identity of currents with $\cL^1_{\loc}$ coefficients.

Hence by (\ref{e54})
$$\left(\frac{1}{\lambda^{k+1}} P_{k+1}^* \omega^n - 
\frac{1}{\lambda^{k}} P_{k}^* \omega^n,\phi\right) =
\frac{1}{\lambda^{k}} \left(P_{k}^* h - P_{k}^* 
\omega^n,\phi\right)\leq \left(\frac{a}{\lambda}\right)^k.$$
for $k$ sufficiently large. Therefore the sequence 
$\{(\frac{1}{\lambda^k} P_k^*\omega^n,\phi)\}$ is Cauchy. Since $\phi$ 
is arbitrary, it follows that $\frac{1}{\lambda^k} P_k^*\omega^n$ 
converges to a measure $\mu.$ %(see    )

Let $\mu_k=\frac{1}{\lambda^k} P_k^*\omega^n.$ Since the $\mu_k$ are 
probability measures, so is $\mu$. The last conclusion of Theorem~\ref{th2} 
then follows from  Theorem~\ref{th1}.\qed

\section*{Acknowledgements} Part of the work on this paper was done while 
both authors were at  the Mathematical 
Sciences Research Institute in Berkeley. This work was started when 
the first author
was visiting Johns Hopkins University. The authors would like 
to thank these organizations for their hospitality and support.


\begin{thebibliography}{MMM}

\bibitem[BLS1]{BLS1} {\sc E. Bedford, M. Lyubich, J. Smillie.}
{\it Polynomial diffeomorphisms of $\C^2$. The measure of
maximal entropy and laminar currents.} Invent.\ Math.\ 112, 1993, 
77--125. 

\bibitem[BLS2]{BLS2} {\sc E. Bedford, M. Lyubich, J. Smillie.}
{\it Distribution of periodic points of polynomial diffeomorphisms of 
$\C^2$.} Invent.\ Math.\ 114, 1993, 277--288. 

\bibitem[BS1]{BedfordSmillie1} {\sc E. Bedford, J. Smillie.}
{\it Polynomial diffeomorphisms of $\C^2:$ currents,
equilibrium measure and hyperbolicity.} Invent.\ Math., 87,
1990, 69--99.

\bibitem[BS2]{BedfordSmillie2} {\sc E. Bedford, J. Smillie.}
{\it Polynomial diffeomorphisms of $\C^2:$ stable manifolds and
recurrence.} J. of Amer.\ Math.\ Soc., 4, 1991, 657--679.

\bibitem[BS3]{BedfordSmillie3}{\sc E. Bedford, J. Smillie.} 
{\it Polynomial diffeomorphisms of $\C^2:$ Ergodicity, exponents and
entropy of the equilibrium measure.} Math.\ Ann.\ 294, 1992,
395--420. 

\bibitem[BT]{BT} {\sc E. Bedford, B. A. Taylor.}
{\it A new capacity for plurisubharmonic functions.} Acta Math.,
149, 1982, 1--40.

\bibitem[Br]{Brolin} {\sc H. Brolin.} {\it Invariant sets under
iteration of rational functions.} Ark.\ Math., 61, 1965, 103--144.

\bibitem[De]{De} {\sc J.-P. Demailly.} {\it Monge-Amp\`ere
operators, Lelong numbers and intersection theory.} Complex Analysis
and Geometry, Univ.\ Series in Math., V. Ancona \& A. Silva, eds., Plenum
Press, New-York, 1993, 115--193.

\bibitem[Fe]{Fe}{\sc H. Federer.} {\it Geometric Measure Theory.}
Springer-Verlag, New York, 1969.

\bibitem[FS1]{FornaessSibony} {\sc J. E. Fornaess, N.  Sibony.} {\it Complex
Henon mappings in $C^2$ and Fatou-Bieberbach domains.} Duke Math.\ J., 65,
1992, 345--380.

\bibitem[FS2]{FS0} {\sc J. E. Fornaess, N.
Sibony.} {\it Complex dynamics in higher dimensions, I.}
Complex analytic methods in dynamical systems (Rio de Janeiro, 1992). 
Ast\'erisque No.\ 222, 1994, 201--231.

\bibitem[FS3]{FS1}{\sc J. E. Fornaess, N. Sibony.} {\it Lectures on
complex dynamics in higher dimensions.}  Notes partially written by 
Estela A. Gavosto. NATO Adv.\ Sci.\ Inst.\ Ser.\ C Math.\ Phys.\ Sci., 
439, Complex potential theory (Montreal, PQ, 1993),
Kluwer Acad.\ Publ., Dordrecht, 1994, 131--186.

\bibitem[FS4]{FS2} {\sc J. E. Fornaess, N. Sibony.}  {\it Complex dynamics in
higher dimensions, II.} Modern methods in complex analysis (Princeton, NJ,
1992), Ann.\ of Math.\ Stud., 137, Princeton Univ.\ Press,  Princeton, NJ,
1995,  135--182.

\bibitem[FS5]{FS5} {\sc J. E. Fornaess, N. Sibony.}  {\it Oka's inequality 
for currents and applications.} Math.\ Ann., 301, 1995,
399--419.

\bibitem[FLM]{FreireLopesMane} {\sc A. Freire, A. Lopes, R.
Ma\~n\'e.} {\it An invariant measure for rational maps.} Bol.\
Soc.\ Bras.\  Mat., 6, 1983, 45--62.

\bibitem[Fr]{Friedland} {\sc S. Friedland.} {\it Entropy of polynomial and
rational maps.} Ann. Math. 133, 1991, 359--368.

\bibitem[Fu]{Fu} {\sc W. Fulton.} {\it Intersection theory.\/}
Springer-Verlag, New York, 1984.

\bibitem[GH]{GH} {\sc P. Griffiths, J. Harris.} {\it Principles of algebraic 
geometry.} Wiley, New York, 1978.

\bibitem[GK]{GK}{\sc P. Griffiths, J. King.} {\it Nevanlinna theory and
holomorphic mappings between algebraic varieties.\/} Acta Math.\ 130 (1973),
145--220.

\bibitem[Ha]{Ha}{\sc R. Hartshorne.} {\it Algebraic geometry.\/}
Springer-Verlag, New York, 1977.

\bibitem[HP]{HubbardPapadopol} {\sc J. H. Hubbard, P. Papadopol.} 
{\it Superattractive fixed points in $\C^n$}. Indiana Univ.\
Math.\ J., 43, 1994, 321--365.

\bibitem[Ly]{Lyubich} {\sc M. Yu. Lyubich.} {\it
Entropy properties of rational endomorphisms of the Riemann
sphere.} Ergodic Theory and Dynamical Systems, 3, 1983, 351--385.

\bibitem[RS]{RS}{\sc A. Russakovskii, M. Sodin} {\it Equidistribution
for sequences of polynomial mappings.} Indiana Univ.\ Math. 
J., 44, 1995, 850--882.

\bibitem[Sh]{Shabat} {\sc B. V. Shabat.} {\it Distribution of
values of holomorphic mappings.} Amer.\ Math.\ Soc., Providence,
R.I., 1985.

\bibitem[Sk1]{Skoda-1} {\sc H. Skoda.} {\it Sous-ensembles
analytiques d'ordre fini ou infini dans $\C^n$.}
Bull.\ Soc.\ Math.\ France, 100, 1972, 353--408.

\bibitem[Sk2]{Skoda-2} {\sc H. Skoda.} {\it Nouvelle m\'ethode pour
l'\'etude des potentiels associ\'es aux ensembles analytiques.}
S\'eminaire P. Lelong (Analyse) 1972 - 1973. Lecture Notes in Math., 410,
Springer-Verlag, New York, 1974, 117--141.

\bibitem[So]{Sodin} {\sc M. L. Sodin.} {\it Value distribution
of sequences of rational mappings.} Advances in Soviet Math.,
11, Entire and Subharmonic Functions, B. Ya.\ Levin, ed., 1992.

\end{thebibliography}
\end{document}